\documentclass[10pt]{amsart}
\usepackage{amsmath,amsthm,amssymb}
\usepackage{graphicx}
\usepackage[margin=1.5in]{geometry}

\begin{document}

%
%

\newcommand{\R} {{\mathbb R}}                              
\newcommand{\Z} {{\mathbb Z}}                              
\newcommand{\N}{{\mathbb N}}                               

\newcommand{\ep} {{\epsilon}}                              

\newcommand{\ao} {{a^o_k}}
\newcommand{\bo} {{b^o_k}}
\newcommand{\au} {{a^u_k}}
\newcommand{\bu} {{b^u_k}}

\newcommand{\Lk}{{\mathcal L}}
\newcommand{\Pa}{{\mathfrak J}}                              
\newcommand{\B}{{\mathcal B}}                              
\newcommand{\stB}{\B^{st}}                                 

\newcommand{\J}{{\mathcal J}}                              
\newcommand{\stab}{{\mathcal S}}                                
\newcommand{\con}{{\mathcal C}}                                 
\newcommand{\Jug} [1] {\langle \ {#1} \ \rangle}           

\newcommand{\C}[1] {C_{#1,n}}
\newcommand{\Cin}[1] {C^{-1}_{#1,n}}
\newcommand{\sig} [1] {\sigma_{#1}^{}}          
\newcommand{\sigin} [1] {\sigma_{#1}^{-1}}      

\newcommand{\st} {\:\: | \:\:}
\newcommand{\ore} {\ \ {\it or} \ \ }
\newcommand{\oand} {\ \ {\it and} \ \ }

%
%

\theoremstyle{plain}
\newtheorem*{prob}{Problem}
\newtheorem{thm}{Theorem}
\newtheorem*{arthm}{Artin's Theorem}
\newtheorem*{mathm}{Markov's Theorem}
\newtheorem*{jrules}{Juggling Rules}
\newtheorem*{crules}{Crossing Rules}
\newtheorem{prop}[thm]{Proposition}
\newtheorem{cor}[thm]{Corollary}
\newtheorem{lem}[thm]{Lemma}
\newtheorem{conj}[thm]{Conjecture}

\theoremstyle{definition}
\newtheorem{defn}[thm]{Definition}
\newtheorem{exmp}[thm]{Example}

\theoremstyle{remark}
\newtheorem*{rem}{Remark}
\newtheorem*{nota}{Notation}
\newtheorem*{ack}{Acknowledgments}
\numberwithin{equation}{section}

\subjclass[2000]{Primary 57M25, Secondary 52C20}

\title{Juggling braids and links}

\author{Satyan L.\ Devadoss}
\address{Satyan L.\ Devadoss: Williams College, Williamstown, MA 01267}
\email{satyan.devadoss@williams.edu}

\author{John Mugno}
\address{John Mugno: University of Maryland, College Park, MD 20742}
\email{jmugno@math.umd.edu}

\begin{abstract}
Using a simplistic model of juggling based on physics, a natural map is constructed from the set of periodic juggling patterns (or \emph{site swaps}) to links.   We then show that all topological links can be juggled.
\end{abstract}

\maketitle

\baselineskip=15pt

%
%
\section{Juggling sequences}

The art of juggling has been around for thousands of years.  Over the past quarter of a century, the interplay between juggling and mathematics has been well studied.  There has even been a book \cite{pol} devoted to this relationship, dealing with several combinatorial ideas.  Numerous juggling software is also available; in particular, Lipson and Wright's elegant and wonderful \emph{JuggleKrazy} \cite{jk} program helped motivate much of this paper.  Most of the information useful to juggling can be accessed via the \emph{Juggling Information Service} webpage~\cite{jis}.  The goal of this paper is to construct and study a  map from juggling sequences to topological braids.  An early form of this idea providing motivation can be found in the work of Tawney \cite{taw}, where he looks at some classic juggling patterns.

In our discussion, we remove everything that is not mathematically relevant. Thus, assume the juggler in question is throwing identical objects, referred to as \emph{balls}.  By convention, there are some basic rules we adhere to in juggling.
\begin{enumerate}
\item[J1.] The balls are thrown to a constant beat, occurring at certain equally-spaced discrete moments in time.
\item[J2.] At a given beat, at most one ball gets caught and then thrown instantly.
\item[J3.] The hands do not move while juggling.
\item[J4.] The pattern in which the balls are thrown is periodic, with no start and no end to this pattern.
\item[J5.] Throws are made with one hand on odd-numbered beats and the other hand on even-numbered beats.
\end{enumerate}
A throw of a ball which takes $k$ beats from being thrown to being caught is called a \emph{$k$-throw}.
Condition J5 implies that when $k$ is even (or odd), a $k$-throw is caught with the same (or opposite) hand from which it was thrown.  Thus, a $5$-throw starting in the left hand would end in the right hand $5$ beats later, while a $4$-throw starting in the left hand would end back in the left hand after $4$ beats.  In this notation, a $0$-throw is a placeholder so that an empty hand can take an action, where no ball gets caught or thrown on that beat.  The following definition from \cite{begw} is used to embody some of the rules above.

\begin{defn}
\label{d:functiondef}
A \emph{juggling pattern} is a bijection
$$f: \Z \rightarrow \Z : t \to t + df(t)$$
where $n \in \N$ and $df(t+n) = df(t) \geq 0$.
A \emph{juggling sequence} (or a \emph{site swap}) is the sequence $\Jug{df(0), \, df(1), \, \cdots, \, df(n-1)}$ arising from a juggling pattern.
\end{defn}

\noindent The number $df(t)$ is the \emph{throw value} at time $t$, and the number $n$ associated to a juggling pattern is its \emph{period}.  Thus, a juggling sequence is simply used to keep track of successive throw values.  Under this terminology, the sequence $\Jug{n}$ yields the $n$-ball \emph{cascade} for $n$ odd, and the $n$-ball \emph{fountain} for $n$ even; the $n$-ball \emph{shower} is given by the $\Jug{1, \, 2n-1}$ pattern.

\begin{figure}[h]
\includegraphics[width=.8\linewidth] {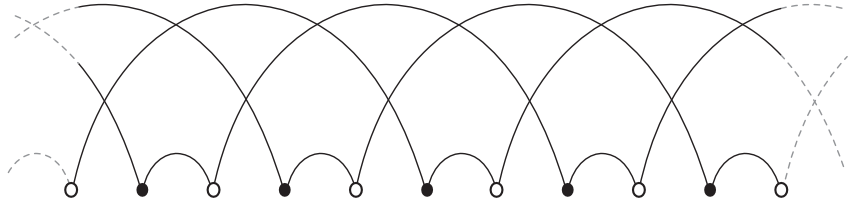}
\caption{Juggling diagram of the $\Jug{5, \, 1}$ sequence.}
\label{f:51diag}
\end{figure}

A juggling sequence can be represented graphically in several ways.  A common way is the \emph{juggling diagram}, with the height of the balls in the vertical direction drawn with respect to time.   Figure~\ref{f:51diag} shows the juggling diagram of the $\Jug{5, \, 1}$ sequence, the lines in this diagram corresponding to the path of the juggled balls traced out over time.
We denote odd and even numbered beats by solid and open circles, distinguishing the two hands.  The juggling pattern associated to the $\Jug{5, \, 1}$ sequence is
\begin{equation*}
f(t) =
\begin{cases}
\ t + 5 & \hspace{10 pt} \mbox{if}  \ t \equiv 0 \mod 2 \\
\ t + 1 & \hspace{10 pt} \mbox{if}  \ t \equiv 1 \mod 2.
\end{cases}
\end{equation*}
It is straightforward to realize that the number of balls used in a juggling pattern $f$ is the number of orbits determined by $f$.  For example, the sequence $\Jug{5, \, 1}$ is a $3$-ball pattern, as seen in Figure~\ref{f:51diag}.
A natural question to ask is which sequence of numbers provide valid juggling sequences. Buhler, Eisenbud, Graham and Wright \cite{begw} provide an elegant criterion.

\begin{thm} \cite{begw}
\label{t:graham}
A sequence $\Jug{h(0), \, h(1), \, \ldots, \, h(n-1)}$ of nonnegative integers is a valid juggling sequence if and only if
\begin{enumerate}
\item the average of the collection $\{ \ h(i) \ \}$ is some integer b, and
\item $\{ \ h(i) + i \mod n \ \}$ is a permutation of \ $\{0, 1, \cdots, n-1\}$.
\end{enumerate}
In this case, the sequence describes a valid $b$-ball juggling sequence.
\end{thm}

%
%
\section{Ladder diagrams}

\begin{figure}[h]
\includegraphics[width=.8\linewidth] {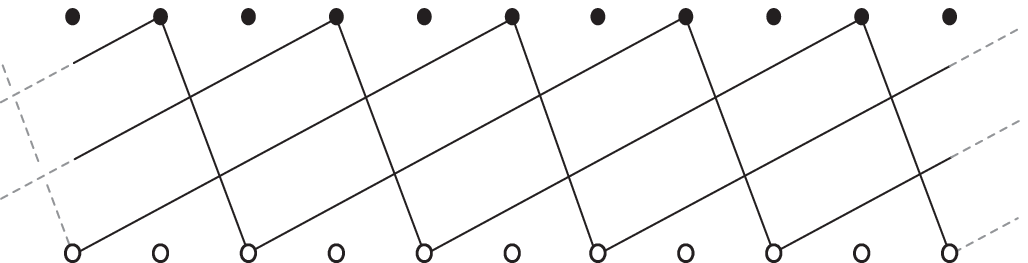}
\caption{Ladder diagram of the $\Jug{5, \, 1}$ sequence.}
\label{f:51ladder}
\end{figure}

With the goal of constructing braids in mind, a slightly different graphical representation is preferable.  In the jugglers terminology, this is called the \emph{ladder diagram}, with the distance between the hands shown in the vertical direction, drawn with respect to time. Figure~\ref{f:51ladder} shows the ladder diagram of the $\Jug{5, \, 1}$ sequence, where one hand is the bottom line and the other is the top.
The viewer is looking at the juggler from above, and the straight lines correspond to the path of the juggled balls traced out over time.  These paths are the parabolic arcs  which appear in the juggling diagram shown in Figure~\ref{f:51diag} above.

\begin{figure}[h]
\includegraphics[width=.8\linewidth] {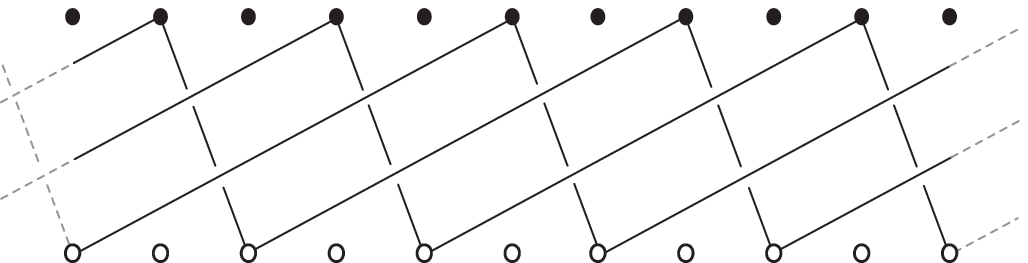}
\caption{Braid diagram of the $\Jug{5, \, 1}$ sequence.}
\label{f:51braid}
\end{figure}

Indeed, from basic physics, the height of a $k$-throw is proportional to $k^2$.   Thus, it is easy to show for a crossing appearing in a ladder diagram, a line representing the trace of ball $x$ will \emph{cross over} a line representing the trace of ball $y$ if and only if $x$ has a higher throw value than $y$.  Figure~\ref{f:51braid} shows the \emph{braid diagram} of the $\Jug{5, \, 1}$ sequence, which encodes crossing information into the ladder diagram.

\begin{rem}
The \emph{physics} of juggling is quite interesting.  We refer the reader to the pioneering work of Magnusson and Tiemann \cite{mt} for details, where dwell times, error margins and angle variations are discussed.
\end{rem}

Before proceeding, two problematic situations need to be addressed, both involving collisions of balls under our simplistic model.  First, interactions between odd throws of the same value need to be discussed.  As mentioned above, basic physics guarantees that odd throws of differing values will not collide in our model.  However, physics also implies two odd throws of the same value which cross in the juggling diagram will collide in reality.  This is because both will have the same throw height, proportional to the square of their throw value (and also from reasons involving symmetry).  Thus, under a simplistic juggling model, the classic juggling sequence $\Jug{3}$ (the $3$-ball cascade) cannot be juggled without the balls colliding in midair.

The second problem involves interactions between even throws.  For example, under our juggling model, a classic juggling sequence such as $\Jug{4}$ (called the \emph{$4$-ball fountain}) cannot be juggled without the balls colliding in midair.  Since a $4$-throw is even, it is thrown vertically in the air, to be returned to the same hand in $4$ beats.  Thus any other even throw from the same hand will cause a collision if it is not thrown far enough apart in time.  This is certainly true of $\Jug{4}$.

Both of these problems can be resolved by adding a bit of realism to our juggling model.  One lesson learned from actual juggling is that the hands move slightly during catches and throws.  Indeed, a catch is made by moving the corresponding hand slightly out of standard position, away from the other hand.  After being caught, the ball is then carried temporarily as it is moved slightly in (closer towards the other hand) before release.\footnote{This is sometimes done in the opposite fashion, with a hand catching on the inside and throwing on the outside, called \emph{reverse juggling}.}
Thus, the following changes are made to the juggling rules, resulting in a modified model.
\begin{enumerate}
\item[J$2^\prime$.] At a given beat, at most one ball gets caught (slightly before the beat) and then thrown (slightly after the beat).
\item[J$3^\prime$.] The hands move slightly such that a ball is caught on the outside and thrown from the inside.
\end{enumerate}

\noindent We show that this minor modification is enough to resolve the collision issues present.

%
%
\section{Dwell and Carry}

In order to quantify our measurements, we introduce \emph{dwell} and \emph{carry} values. The dwell value $d$ is the amount of time (measured with respect to the beat of the throws) a ball is held between catching it and throwing it.  The carry value $c$ is the distance each hand moves from catching to throwing.  For our purposes, let $c > 0$ and $d > 0$, both being very small values of our choosing.

\begin{figure}[h]
\includegraphics [width=.9\linewidth]{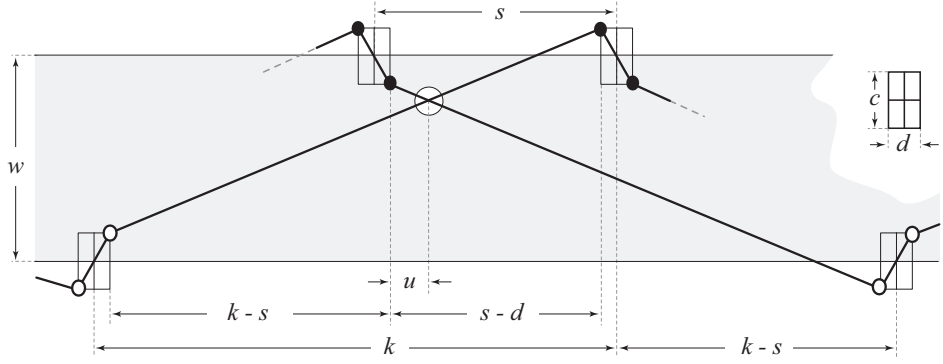}
\caption{Details of carry and dwell.}
\label{f:proof}
\end{figure}

Figure~\ref{f:proof} provides the ladder diagram of the basic set-up:  Two odd $k$-throws are examined with dwell and carry values.  The difference of beats between the second ball thrown and the first ball caught is denoted by $s$, with an actual difference of $s-d$ taking the dwell value into consideration.  Let $u$ be the difference in time between the throw of the second ball and the intersection of the two paths from the perspective of the ladder diagram.  We can calculate the actual value of $u$ by solving for the intersection of the two lines, obtaining
$$u = \frac{1}{2}(s- c^* - d),$$
where
$$c^* = \frac{c}{w}(k-d).$$
We need to find out which ball will be higher at this point of intersection.

Since the height is parabolic, with maximum height reached at the midpoint of the throw distance, we need to compute which throw is closest to the midpoint value at the intersection.  The second throw has horizontal distance $u$.  The first throw has horizontal distance $k-s+u$ to the point of intersection.  From symmetry, this is the same height as a distance of
$$(k-d) - (k-s+u) = s-d-u = \frac{1}{2}(s+c^*-d).$$
The midway distance of maximal height is $\frac{1}{2}(k-d)$.  Since we are assuming $k > s$ and since $k$ and $s$ are natural numbers, then by the Archimedean principle we can choose a value for $c$ such that $k > s + c^*$.  Thus,
$$\frac{1}{2}(k-d) > \frac{1}{2}(s+c^*-d) > \frac{1}{2}(s-c^*-d) = u.$$
In other words, the first (earlier) throw will have a higher height than second (later) one at the crossing.  

\begin{rem}
Notice the arguments above carry through when $d = 0$.  However, for the sake of realism, we choose a small positive value for $d$.
\end{rem}

\begin{thm}
\label{t:collide}
For a given juggling sequence, the balls will not collide for a small enough choice of carry value $c>0$.
\end{thm}

\begin{proof}
Any crossing in the ladder diagram will be of two kinds:  one between two throws (Figure~\ref{f:oddeven}A) or one between a throw and a carry (Figure~\ref{f:oddeven}B).  Any crossing of the latter kind is easily resolved, since  the line traced out by a carry (the movement of a hand carrying a ball) has height zero.  Thus, only crossings between two throws $\alpha$ and $\beta$ need to be addressed.

\begin{figure}[h]
\includegraphics [width=.9\linewidth]{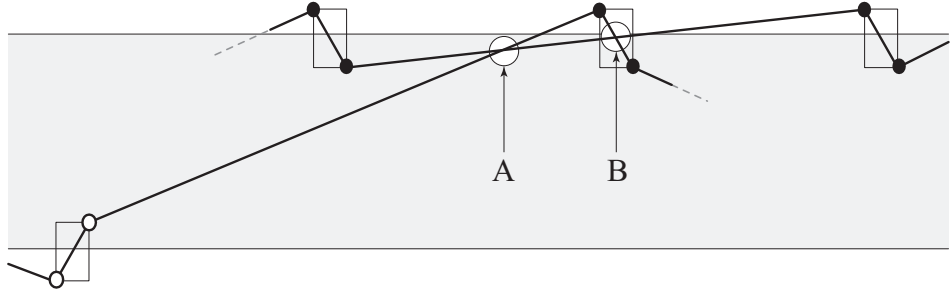}
\caption{Types of intersections.}
\label{f:oddeven}
\end{figure}

When $\alpha > \beta$, the $\alpha$-throw will cross above the $\beta$-throw since the balls will have different heights proportional to the square of their throw values.  When $\alpha = \beta$ is odd,  the argument outlined above (for small enough $c>0$ value) will guarantee that the earlier of the two throws will cross over the latter.  If $\alpha=\beta$ is even, the throws will never cross since they travel along parallel trajectories.

The last possibility is when one value is odd and the other even, as shown in Figure~\ref{f:oddeven}A.  Here, as $c \to 0$, this intersection will occur arbitrarily close to when the odd throw is caught (at height zero) while the even throw is still in the air.  Thus for a small enough value of $c >0$, the balls will not collide.
\end{proof}

\begin{figure}[h]
\includegraphics [width=.8\linewidth]{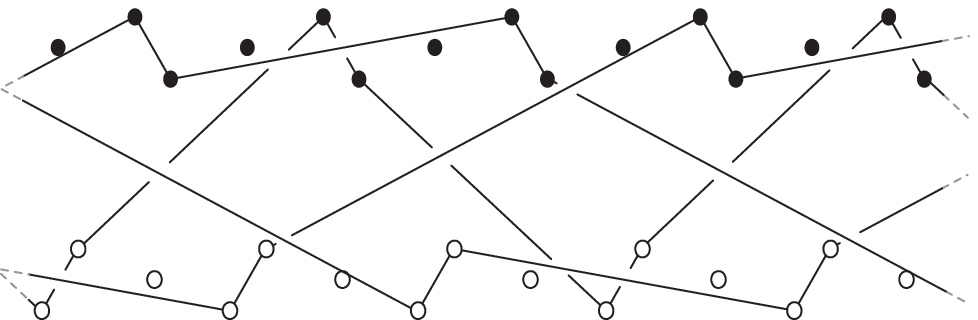}
\caption{Braid diagram of the $\Jug{3, \, 4, \, 5}$ sequence.}
\label{f:345braid}
\end{figure}

Indeed, all crossing issues can be resolved due to the theorem above;  Figure~\ref{f:345braid} shows an example, the braid diagram of the $\Jug{3, \, 4, \, 5}$ sequence.

\begin{rem}
The rest of the paper will mostly be concerned with crossings of odd throws of differing values.  Thus, although dwells and carries can occur in our juggling model, we will henceforth visually represent the braid diagrams without them.
\end{rem}

%
%
\section{Braids and Links}

Braids have a rich history in mathematics, appearing in numerous areas.  We refer the reader to Adams \cite[Chapter 5]{ada} for an elementary introduction to this subject.  Roughly, an \emph{$n$-braid} is a set of $n$ disjoint arcs running between two vertical bars where every vertical plane between the two bars must intersect each arc exactly once.  The endpoints of the arcs are given an ordered labeling.

Starting with a trivial $n$-braid, let $\sig{i}$ denote the braid where the strand at the $i$-th position crosses \emph{over} the strand at the $(i+1)$-st position.  Similarly, let $\sigin{i}$ represent the braid with the crossing of the strand at the $i$-th position \emph{under} the $(i+1)$-st position.  Since any braid can be obtained by repeatedly crossing adjacent strands, every braid can be expressed as a \emph{word}, an ordered collection of $\sig{}$ elements.

The $n$-braids have an algebraic group structure, where the composition map of two braids $w_1 \cdot w_2$ simply attaches the $n$ ordered endpoints of the arcs of $w_1$ to the $n$ starting points of the arcs of $w_2$.   We denote the group of $n$-braids by $\B_n$, calling it the \emph{braid group}.  In actuality, elements of the braid group are \emph{equivalence classes} of braids.  The following theorem by Artin (1926) gives the relations needed to identify isomorphic braids:

\begin{arthm}
The braid group $\B_n$ is generated by $\{\sig{1}, \cdots, \sig{n-1}\}$ with the relations:
\begin{enumerate}
\item[R1.] $\sig{i} \cdot \sig{i+1} \cdot \sig{i} = \sig{i+1} \cdot \sig{i} \cdot \sig{i+1}$.
\item[R2.] $\sig{i} \cdot \sig{j} = \sig{j} \cdot \sig{i}$ \ \ if \ $|i-j|>1$.
\end{enumerate}
\end{arthm}

A natural operation to perform on braids is their \emph{closure}, identifying the endpoints of the strands of the braids.  The result after closing a braid is a \emph{solid torus braid}, a braid in a solid torus where each meridional disk intersects each strand of the braid exactly once.
Removing the restriction of being in a solid torus allows a closed braid to become a \emph{link} in $\R^3$.
Markov (1935) introduced two new relations which algebraically encapsulate our discussion. Let $w$ be a word in $\B_n$.
\begin{enumerate}
\item[R3.] (Conjugation) $w = \sig{i} \cdot w \cdot \sigin{i}$.
\item[R4.] (Stabilization) $w = w \cdot \sig{n}$ \ \ where $\sig{n} \in \B_{n+1}$.
\end{enumerate}

\begin{mathm}
The braid group up to the conjugation relation describes the solid torus braids.  The braid group up to both conjugation and stabilization describes links.
\end{mathm}

\noindent Thus, the set of all solid torus braids (and links) can be generated by $\sig{}$ elements and their inverses.

The smallest part of the juggling diagram which can be used to tile it is called the \emph{fundamental chamber}.  More precisely, a juggling sequence has \emph{exact period} $n$ if it does not have another period $m$ for a divisor $m$ of $n$.  Thus, a fundamental chamber for a juggling sequence with exact period $n$ has length $2n$ when $n$ is odd and $n$ when it is even.  This asymmetry between odd and even length patterns arises from throwing from right and left hands.  After an odd number of throws, the next throw will be from the opposite hand; thus, doubling the pattern will be the first repeatable set of throws that ends on the correct hand.

\begin{defn}
Let $\J$ be the map from juggling sequences to solid torus braids, defined by taking the closure of a fundamental chamber of the braid diagram.
\end{defn}

\begin{rem}
Theorem~\ref{t:collide} guarantees $\J$ to be well-defined.  Note that there is no canonical map from juggling sequences to braids.   The translation of the fundamental chamber produces alterations of the starting and ending points of the braid.  Indeed, taking the closure of the chamber removes all ambiguity, resulting in solid torus braids.
\end{rem}

\begin{figure}[h]
\includegraphics[width=\linewidth] {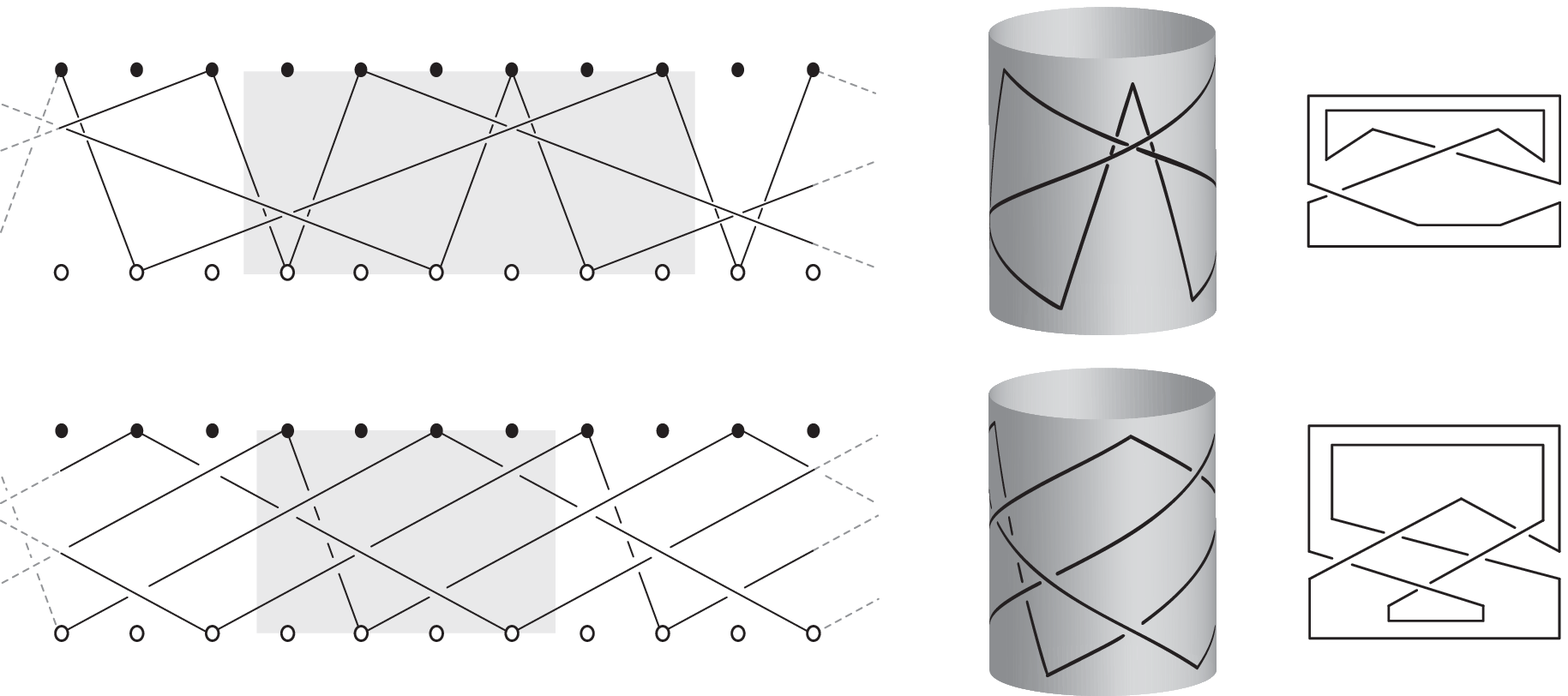}
\caption{The juggling sequences $\Jug{7, \, 1, \, 1}$ and $\Jug{5, \, 5, \, 5, \, 1}$.}
\label{f:knots}
\end{figure}

Figure~\ref{f:knots} shows examples of $\J$ evaluated on  sequences $\Jug{7, \, 1, \, 1}$ and $\Jug{5, \, 5, \, 5, \, 1}$, respectively.  A fundamental chamber in each braid diagram is shaded. By placing the solid torus braids in $\R^3$ (allowing stabilization), then $\Jug{7, \, 1, \, 1}$ maps to the unknot and $\Jug{5, \, 5, \, 5, \, 1}$ maps to the trefoil.\footnote{This could be the right or the left trefoil, depending on which hand releases the ``1'' throw.  In other words, the map $\J$ is up to reflections of links.}

The natural question is to ask which links arise from juggling patterns.  The following is the key result, the proof of which will consume most of the remaining paper.

\begin{thm}
\label{t:bigthm}
The map $\J$ is surjective.
\end{thm}

\noindent This, along with Markov's theorem, results in the following:

\begin{cor}
For every link, there exists a juggling sequence which maps to the link.
\end{cor}

%
%
\section{Constructing a trivial braid}

Our method of attacking Theorem~\ref{t:bigthm} is by constructing a juggling pattern for every solid torus braid.  We begin with the trivial solid torus $n$-braid, which will serve as a building block.  Before providing the juggling pattern, it is first described from the perspective of the balls (or strands of the braid).  The first ball of the $n$ balls is thrown at a constant beat of $3$ units, being thrown at $t=0$.  The second ball is thrown at a constant beat of $3^3$ units, being thrown at $t=1$.  Let $\alpha_1 = 0$ and let $\alpha_k = 3^0 + 3^2 + \cdots + 3^{2(k-2)}$ for $k \geq 2$.  In general, the $k$-th ball is thrown at a constant beat of $3^{2k-1}$, being thrown at time $t=\alpha_k$.  We denote this pattern as $I_n$, where $n$ corresponds to the number of balls thrown.

\begin{figure}[h]
\includegraphics[width=.8\linewidth]{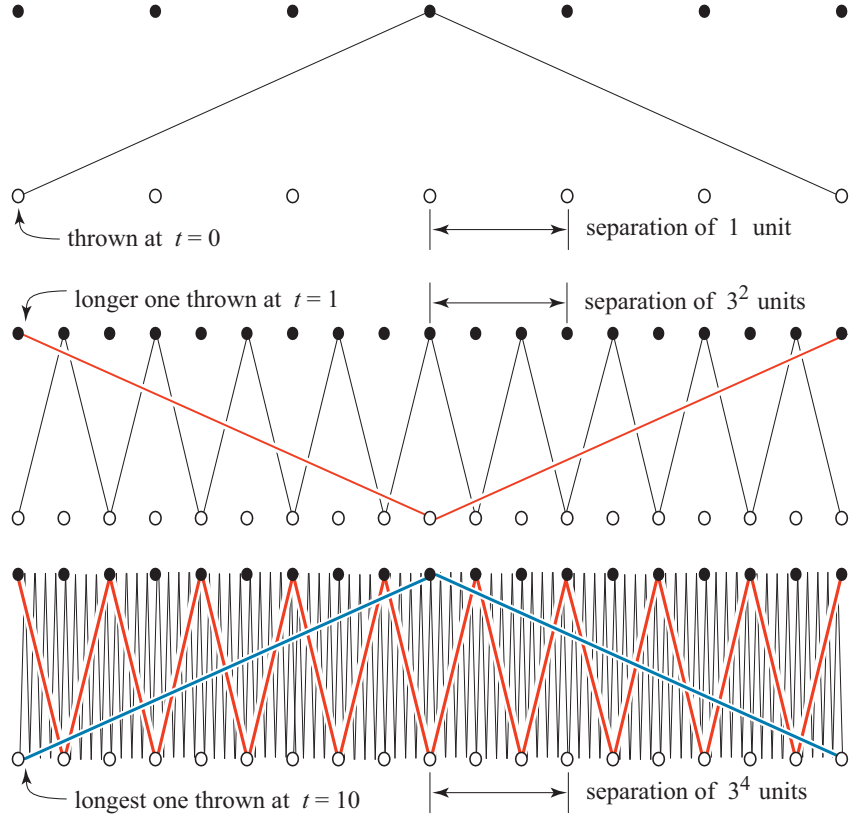}
\caption{The braid diagrams for $I_1$, $I_2$ and $I_3$.}
\label{f:In}
\end{figure}

Figure~\ref{f:In} shows the braid diagrams for $I_1$, $I_2$ and $I_3$ respectively.  We refer to the strand traced out by the $k$-th ball as the \emph{$k$-strand}.
Notice how a given strand always crosses over or always crosses under any another strand.
The juggling pattern can then be defined as
\begin{equation}
\label{e:iden}
I_n(t) =
\begin{cases}
\ t + 3^{2k-1} & \hspace{10 pt} \mbox{if}  \ t \equiv (\alpha_k + m \cdot 3^{2k-1}) \mod 3^{2n-1} \\
\ t + 0 & \hspace{10 pt}  \mbox{otherwise}
\end{cases}
\end{equation}
where $k \in \{1, 2, \ldots, n\}$ and $m \in \Z$.
Here, $k$ keeps track of the balls and $\alpha_k + m \cdot 3^{2k-1}$ the beats in which the $k$-th ball is thrown.

\begin{thm}
The function $I_n(t)$ is a valid juggling pattern.
\end{thm}

\begin{proof}
We need to show that Eq.~\eqref{e:iden} satisfies the two criteria set forth in Theorem~\ref{t:graham}.  The set of beats in which ball $k$ is thrown (with throw length $3^{2k-1}$) is
\begin{equation}
\label{e:throws}
\{\ \alpha_k + m \cdot 3^{2k-1} \mod 3^{2n-1} \ \},
\end{equation}
where $m \in \Z$.  When $m = 3^{2(n-k)}$, then $\alpha_k + m \cdot 3^{2k-1} \equiv \alpha_k \mod 3^{2n-1}$ and this is clearly the smallest such positive value of $m$ for which this property holds.
Thus, there are $3^{2(n-k)}$ throws of length $3^{2k-1}$.  Taking the average of all the throw lengths yields
\begin{equation}
\label{e:ave}
\frac{1}{3^{2n-1}} \sum_{k=1}^n 3^{2(n-k)} \cdot 3^{2k-1} = \frac{1}{3^{2n-1}} \sum_{k=1}^n 3^{2n-1} = n,
\end{equation}
meeting the first criteria.  So the number of balls in this pattern is $n$, as desired.

For the second condition, we need to show
\begin{equation}
\label{e:balance}
\{\ dI_n(i) + i \mod 3^{2n-1} \} = \{\ i \mod 3^{2n-1} \},
\end{equation}
where $i \in \Z$ marks the position and $dI_n(i)$ the throw value.  This equality clearly holds at positions when no balls are thrown, since the throw value is $0$.  Look at each ball separately and recall that the set of all throw positions of ball $k$ is given by \eqref{e:throws} above, corresponding to the right hand side of Eq.~\eqref{e:balance}.  The left hand side of Eq.~\eqref{e:balance} for a given ball $k$ is defined by
$$\{ \ 3^{2k-1} + (\alpha_k + m \cdot 3^{2k-1})  \mod 3^{2n-1} \ \} \ = \ \{ \ \alpha_k + (m+1) \cdot 3^{2k-1}  \mod 3^{2n-1} \ \}.$$
Since $m \in \Z$, this is a rephrasing of \eqref{e:throws}.
\end{proof}

\begin{thm}
\label{t:trivial}
The image of $I_n$ under $\J$ is the trivial solid torus $n$-braid.
\end{thm}

\begin{proof}
From Eq.~\eqref{e:ave}, $I_n$ is an $n$-ball pattern, so $\J(I_n)$ is a solid torus $n$-braid.  Since all throws of a given ball $k$ are of the same height $3^{2k-1}$, from Theorem~\ref{t:collide}, the $k$-strand will cross over all $j$-strands, when $j < k$ and cross under all $j$-strands when $j > k$. This results in the trivial solid torus braid.
\end{proof}

%
%
\section{Constructing the generators}

There is a well-known method of constructing new juggling patterns from old ones.  Let $\Jug{h(0), \, h(1), \, \cdots, \, h(n-1)}$ be a juggling sequence.  Let $a$ and $b$ be integers such that $0 \leq a < b \leq n-1$ and $b-a \leq h(a)$.  We construct a new sequence which coincides with $h(i)$ on all beats except at $a$ and $b$.  Thus, define $g(i) = h(i)$ for all $i \neq a, b$, and let $g(a) = h(b) + (b - a)$ and $g(b) = h(a) - (b - a)$.

\begin{thm}
\label{t:swap}
$\Jug{g(0), \, g(1), \, \cdots, \, g(n-1)}$ is a juggling sequence.
\end{thm}

\begin{proof}
The two conditions of Theorem~\ref{t:graham} need to be satisfied.  Clearly the average of the collection $\{ \ g(i) \ \}$ is the same as the average of $\{ \ h(i) \ \}$, since $g(a)$ adds $b-a$ to its value and $g(b)$ subtracts it.
Moreover $\{ \ g(i) + i \mod n \ \}$ is a permutation of \ $\{0, 1, \cdots, n-1\}$ since
$g(a) + a = h(b) + b$ and $g(b) + b = h(a) + a$.
\end{proof}

\noindent In other words, the \emph{sites} $h(a)$ and $h(b)$ in which the balls land will \emph{swap} positions.  Indeed, the site swap terminology derives its name from this property.

\begin{figure}[h]
\includegraphics [width=.9\linewidth]{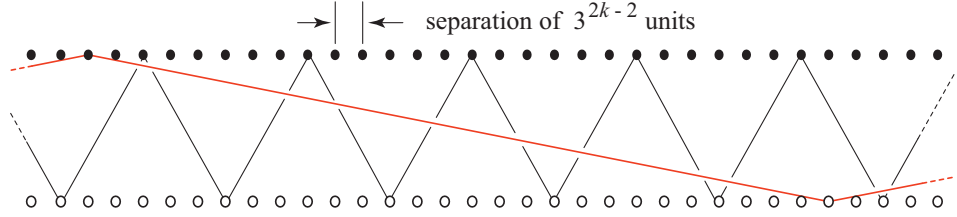}
\caption{Braid diagram only of the $k$ and $k+1$ strands of $I_n(t)$.}
\label{f:trivial}
\end{figure}

In our context, when attempting to construct the generators of the braid group, the creation of a crossing is simply an extension of this idea of swapping sites.  Starting with the identity braid $I_n$, the strands to be crossed are chosen and two of their throws are manipulated (swapped), creating a crossing in an otherwise trivial braid.  Figure~\ref{f:trivial} shows a part of the braid diagram of $I_n(t)$, where only the $k$ and $k+1$ strands are depicted.  The length of beats between two hands is scaled to a factor $3^{2k-2}$.

\noindent {\emph{Over-crossing:}} \newline
Let $\ao = \alpha_{k+1}$ and $\bo = \alpha_k + 3^{2k+1}$.  Then $\bo-\ao = (3^3-1) \cdot 3^{2k-2} < 3^{2(k+1)-1} = dI_n(\ao)$, satisfying conditions of Theorem~\ref{t:swap} above.  We take a fundamental chamber of $I_n(t)$ and swap the two sites $dI_n(\ao)$ and $dI_n(\bo)$, resulting in
\begin{equation*}
\C{k}(t) =
\begin{cases}
\ I_n(\bo) + (\bo-\ao) & \hspace{10 pt} \mbox{if}  \ t \equiv \ao \mod 2 \cdot 3^{2n-1} \\
\ I_n(\ao) - (\bo-\ao) & \hspace{10 pt} \mbox{if}  \ t \equiv \bo \mod 2 \cdot 3^{2n-1} \\
\ I_n(t) & \hspace{10 pt}  \mbox{otherwise.}
\end{cases}
\end{equation*}
Notice the length of its fundamental chamber is $2 \cdot 3^{2n-1}$.

\begin{lem}
$\J(\C{k}(t))$ is the solid torus $n$-braid with its $k$-strand crossing over its $(k+1)$-strand.
\end{lem}

\begin{figure}[h]
\includegraphics [width=.9\linewidth]{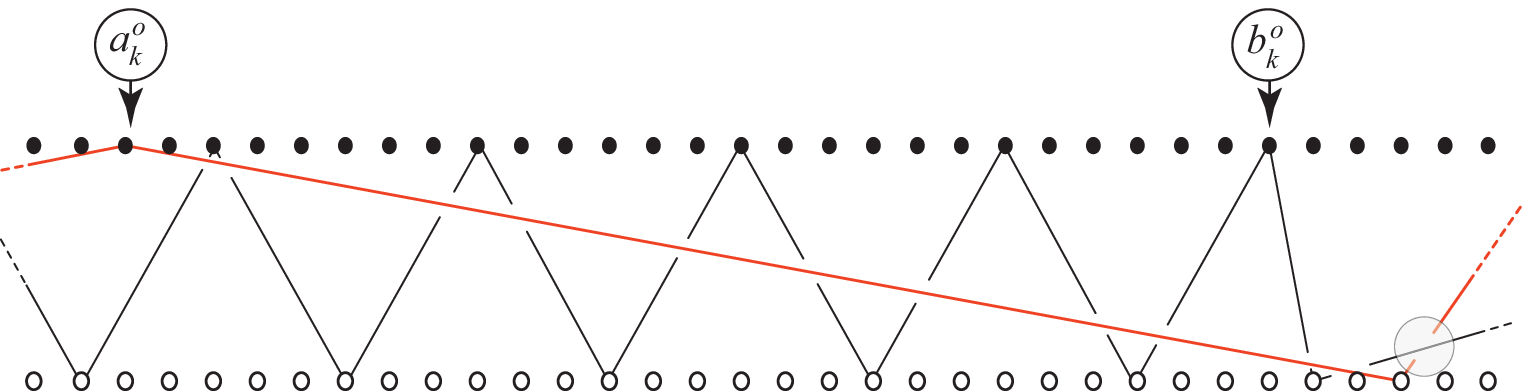}
\caption{Braid diagram only of the $k$ and $k+1$ strands of $\C{k}(t)$.}
\label{f:over}
\end{figure}

\begin{proof}
Theorems~\ref{t:trivial} and \ref{t:swap} guarantee that $\C{k}(t)$ is a juggling pattern.  We need to show the crossing information is as claimed.  The throw value of ball $k+1$ at $\ao$ increases from $3^{2(k+1)-1}$ to
$$3^{2k-1} + (\bo-\ao) = (3^3 +3 -1) \cdot 3^{2k-2}.$$
Crossings will not interfere with other $j$-strands, where $j > k+1$ since
$$3^{2j-1} \geq 3^{2(k+2)-1} = 3^5 \cdot 3^{2k-2} > (3^3 +3 -1) \cdot 3^{2k-2}.$$
Similarly, the throw value of ball $k$ starting at position $\bo$ decreases from $3^{2k-1}$ to
$$3^{2(k+1)-1} - (\bo-\ao) = 3^{2k-2}.$$
So crossings will not interfere with other $j$-strands, where $j < k$.
Thus, the only crossings that can occur are between the $k$ and $k+1$ strands around the swap location. Figure~\ref{f:over} depicts the details:  The $k$-strand at position $\bo + 3^{2k-2}$ has a throw value of $3^{2k+1}$, thus crossing over the $(k+1)$-strand thrown from $\bo + 3^{2k-1}$ with a value of $3^{2k-1}$.
\end{proof}

\noindent {\emph{Under-crossing:}} \newline
This is identical to the situation above, except for different choices of sites to swap.
Let $\au = \alpha_{k+1}$ and $\bu = \alpha_k + 3^{2k-1}$.  Then $\bu-\au = 2 \cdot 3^{2k-2} < 3^{2(k+1)-1} = dI_n(\au)$, satisfying conditions of Theorem~\ref{t:swap} above.   We take a fundamental chamber of $I_n(t)$ and swap the two sites $dI_n(\au)$ and $dI_n(\bu)$, resulting in
\begin{equation*}
\Cin{k}(t) =
\begin{cases}
\ I_n(\bu) + (\bu-\au) & \hspace{10 pt} \mbox{if}  \ t \equiv \au \mod 2 \cdot 3^{2n-1} \\
\ I_n(\au) - (\bu-\au) & \hspace{10 pt} \mbox{if}  \ t \equiv \bu \mod 2 \cdot 3^{2n-1} \\
\ I_n(t) & \hspace{10 pt}  \mbox{otherwise.}
\end{cases}
\end{equation*}

\begin{lem}
$\J(\Cin{k}(t))$ is the solid torus $n$-braid with its $k$-strand crossing under its $(k+1)$-strand.
\end{lem}

\begin{figure}[h]
\includegraphics [width=.9\linewidth]{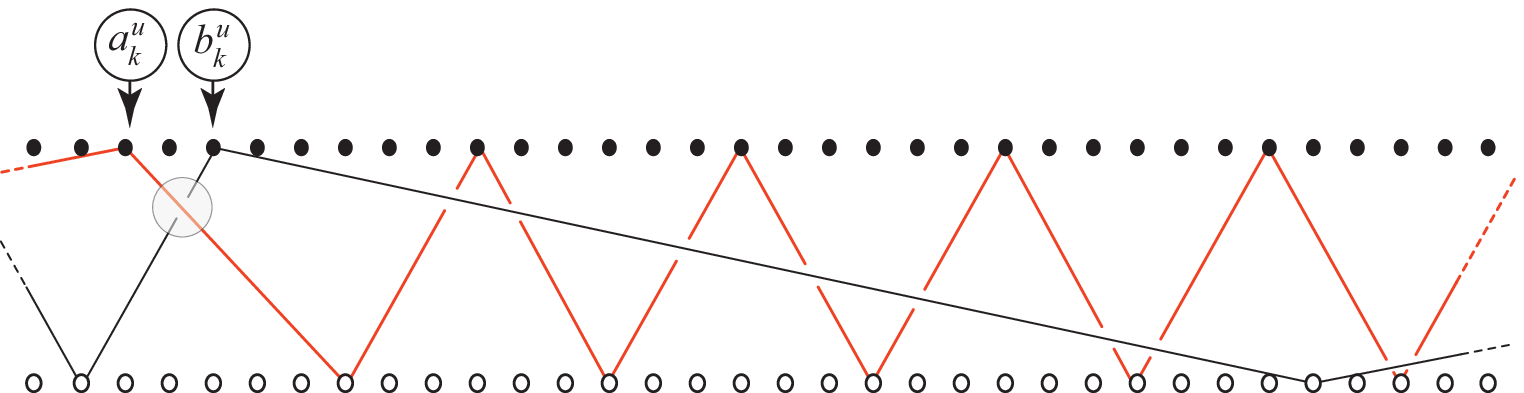}
\caption{Braid diagram only of the $k$ and $k+1$ strands of $\Cin{k}(t)$.}
\label{f:under}
\end{figure}

\begin{proof}
Again, Theorems~\ref{t:trivial} and \ref{t:swap} guarantee $\C{k}(t)$ to be a juggling pattern. The throw value of ball $k+1$ at $\au$ decreases from $3^{2(k+1)-1}$ to
$$3^{2k-1} + (\bu-\au) = 3^{2k-1} + 2 \cdot 3^{2k-2}.$$
Thus crossings will not interfere with other $j$-strands, where $j \leq k$.
Similarly, the throw value of ball $k$ starting at position $\bu$ increases from $3^{2k-1}$ to
$$3^{2(k+1)-1} - (\bu-\au) = 3^{2k+1} - 2 \cdot 3^{2k-2}.$$
So crossings will not interfere with other $j$-strands, where $j \geq k+1$.
Thus, the only crossings that can occur are between the $k$ and $k+1$ strands around the swap location; see Figure~\ref{f:under} for details.  Notice the $k$-strand at position $\au - 3^{2k-2}$ with throw value of $3^{2k-1}$ crosses under the $(k+1)$-strand at position $\au$ having a throw value of $3^{2k-1} + 2 \cdot 3^{2k-2}$.
\end{proof}

To finish the proof of Theorem~\ref{t:bigthm}, we construct a juggling pattern for every solid torus braid.  Let $w=\sigma_{i_1}^* \cdot \sigma_{i_2}^* \cdot \cdots \cdot \sigma_{i_r}^*$ be a word describing a solid torus $n$-braid.  Let
\begin{equation*}
I^w_n(t) =
\begin{cases}
I_n(b^*_{i_j}) + (b^*_{i_j} - a^*_{i_j}) & \hspace{10 pt} \mbox{if}
\ t \equiv a^*_{i_j} + 2 (j-1) \cdot 3^{2n-1} \mod (2r \cdot 3^{2n-1}) \\
I_n(a^*_{i_j}) - (b^*_{i_j} - a^*_{i_j}) & \hspace{10 pt} \mbox{if}
\ t \equiv b^*_{i_j} + 2 (j-1) \cdot 3^{2n-1} \mod (2r \cdot 3^{2n-1}) \\
I_n(t) & \hspace{10 pt}  \mbox{otherwise.}
\end{cases}
\end{equation*}
If $\sigma_i^* = \sig{i}$, then $b^*_i = b^o_i$ and $a^*_i = a^o_i$; similarly, if $\sigma_i^* = \sigin{i}$, then $b^*_i = b^u_i$ and $a^*_i = a^u_i$.  Theorem~\ref{t:trivial} along with the lemmas above guarantee $I^w_n(t)$ to be a juggling pattern; we leave it to the reader to provide details.  We claim that $\J(I^w_n)$ maps to $w$.  In $I^w_n(t)$, $r$ copies of the fundamental chambers of $I_n(t)$ are used, one for each element in $w$; the length of its fundamental chamber is $2r \cdot 3^{2n-1}$.  Each copy is altered by an appropriate swapping of sites corresponding to the generating element in $w$.  This alteration provides the appropriate crossing needed, and provides the proof of Theorem~\ref{t:bigthm}.

\begin{rem}
Choosing the starting position of the fundamental chamber of $I^w_n(t)$ yields an ordering of the \emph{position} of the strands, determining $w$.  This ordering is \emph{not necessarily} the ordering of the $n$ balls, which are used to label the strands.  This is crucial to note, since a site swap switches the \emph{strand} formed by the $k$ ball with the \emph{strand} formed by the $k+1$ ball.  However, since our elements are solid torus braids, we have the ability to slide our strands in the solid torus to the appropriate order (due to conjugation by Markov's theorem).
\end{rem}


%
%
\section{Looking Forward}


Although our construction allows us to prove all links can be juggled, it is far from realistic or efficient.
The throw values were chosen in powers of 3 (a prime) in order to make transparent the construction of a valid juggling sequence, along with isolating crossings of two strands.  This requires a juggling sequence with $n$ balls to have throw values up to $3^{2n-1}$.  Since juggling values of $9$ are near impossible to perform, the method above is certainly not realistic.  We now look at how realism can be introduced and measured.

\begin{figure}[h]
\includegraphics[width=\linewidth]{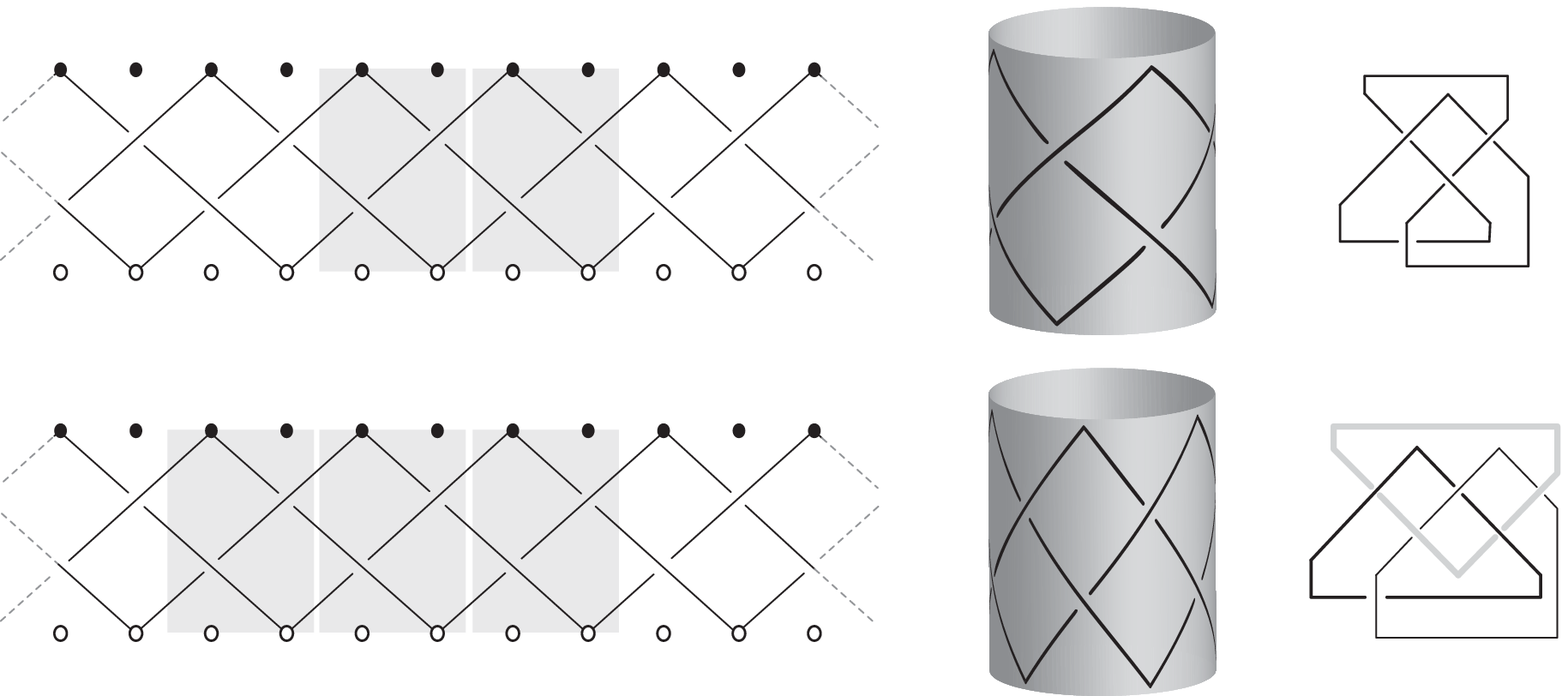}
\caption{$\J^2 \Jug{3}$ and $\J^3 \Jug{3}$ mapping to the figure-eight knot and the Borromean rings.}
\label{f:knot3-2}
\end{figure}

We begin with the map $\J$ from juggling sequences to solid torus braids.  This map is based on taking the closure of the fundamental chamber.  Thus, the classic $\Jug{3}$ cascade sequence, having fundamental chamber of length two, is allowed only to be ``active'' for two beats, resulting in the unknot.  But what if our juggler wishes to juggle longer, for more beats?   Since we want to close the resulting braid, juggling multiple copies of the fundamental chamber can be allowed.  Define $\J^k (f)$ to be the closure of $k$ adjacent copies of the fundamental chamber associated to the juggling pattern $f$.  Figure~\ref{f:knot3-2} shows $\J^2 \Jug{3}$ and $\J^3 \Jug{3}$; the first maps to the figure-eight knot,  the latter to the Borromean rings.

Given a link, we ask to find the best juggling pattern which maps to it.  One can try to measure what ``best'' means by looking at a few factors.  Given a link $l$, let $\Pa(l)$ be the set of juggling patterns $f$ where $\J^k(f)$ maps\footnote{Although $\J^k$ maps a juggling sequence to solid torus braids, we abuse terminology and sometimes refer to the composition of this map with stabilization.} to $l$ for some $k \in \N$.

\begin{defn}
The \emph{ball index}\footnote{Clearly, the \emph{braid index} of a link is less than or equal to the ball index.}   of a link $l$ is the minimum number of balls needed for juggling pattern $f$, for all $f$ in $\Pa(l)$. The \emph{throw index} of a link $l$ is the minimum over all the maximum throw values of a juggling pattern $f$, for all $f$ in $\Pa(l)$.
\end{defn}

It is straight-forward to show that for non-trivial links, the ball index must be greater than one and the throw index must be greater than two.  Consider the trefoil as an example.  A few possible ways to construct it are by $\J^2 \Jug{3}$, $\J^3 \Jug{4, \, 0}$ and $\J \Jug{5, \, 5, \, 5, \, 1}$, as shown in Figure~\ref{f:knots}.  Since $\Jug{4, \, 0}$ is a $2$-ball pattern, the ball index of the trefoil must be two.  Moreover,  $\Jug{3}$ guarantees a throw index of three for the trefoil.

\begin{prob}
Study the properties of the ball index and throw index of links.
\end{prob}

An underlying issue to this problem is understanding $\J$.  We have only been able to prove that $\J$ is surjective, not knowing much more about the map itself.  Based on the figures above, \emph{up to conjugation}, $\J$ maps
\begin{eqnarray*}
\Jug{7, \, 1, \, 1} & \longrightarrow  & \sig{2}\sig{1}\sigin{2}\sigin{1}\sigin{2}\cdot\sigin{1}\sigin{2}\sig{1}\sig{2}\sig{1} \\
\Jug{5, \, 5, \, 5, \, 1} & \longrightarrow & \sig{2}\sig{3}\sig{2}\sig{1}\sigin{2}\sigin{1}\sig{3}  \\
\Jug{3, \, 4, \, 5} & \longrightarrow & \sigin{3}\sigin{1}\sig{3}\sig{2}\cdot\sig{1}\sig{3}\sigin{1}\sigin{2}  \\
\Jug{4} & \longrightarrow & \sigin{1}\sig{3}.
\end{eqnarray*}
The first word is an element in the solid torus braid with three strands; the last three have four strands.
Thus, an algorithmic description of this map is needed, not based on crossing information in figures but on values of the juggling sequences.

\begin{prob}
For any juggling sequence, find an algebraic description of its solid torus braid.
\end{prob}

The version of a juggling pattern used here is sometimes referred to as \emph{vanilla site swap}.  Allowing multiple jugglers, simultaneous throws, and varying throws (reverse juggling) can be encoded into more complicated models; Polster \cite[Chapter 4]{pol} provides a nice introduction to these ideas.  We end with the following:

\begin{prob}
Study braids and links obtained in more complicated juggling models.
\end{prob}

\begin{ack}
We  thank Colin Adams for conversations in the early stages of this work, and Jake Tawney for motivating the problem.  The first author is grateful to the NSF for partially supporting this project with CARGO DMS-0310354.
\end{ack}

%
%
\bibliographystyle{amsplain}

\begin{thebibliography}{XX}
\baselineskip=12pt

\bibitem[1]{ada} C.\ Adams, \emph{The Knot Book}, W.H.\ Freeman, New York, 1994.

\bibitem[2]{begw} J.\ Buhler, D.\ Eisenbud, R.\ Graham, C.\ Wright, Juggling drops and descents, \emph{The American mathematical Monthly} {\bf 101} (1994) 507-519.

\bibitem[3] {jis} Juggling Information Service (JIS), \verb+http://www.juggling.org+

\bibitem[4]{jk} A.\ Lipson and C.\ Wright, \emph{JuggleKrazy} software.

\bibitem[5] {mt} B.\ Magnusson and B.\ Tiemann, The physics of juggling, \emph{Physics Teacher} {\bf 27} (1989) 584-589.

\bibitem[6] {pol} B.\ Polster, \emph{The Mathematics of Juggling}, Springer-Verlag, New York, 2003.

\bibitem[7] {taw} J.\ Tawney, Jugglinks,  Master's Thesis, The Ohio State University, May 2001.
\end{thebibliography}

\end{document}